\documentclass[3p,times]{elsarticle}
\usepackage{amsfonts}
\usepackage{amssymb}
\usepackage{amsmath}
\usepackage{fancyvrb}
\usepackage{graphics}
\usepackage{graphicx}
\usepackage{epsfig}
\usepackage{epstopdf}
\usepackage{lipsum}
\usepackage{subfigure}
\usepackage{leftidx}
\usepackage{multirow}
\usepackage{charter}
\usepackage[charter]{mathdesign}

\usepackage{amssymb}

\begin{document}
\newtheorem{thm}{Theorem}
\newtheorem{lem}[thm]{Lemma}
\newtheorem{defn}[thm]{Definition}
\begin{frontmatter}



\title{A novel second order finite difference discrete scheme for fractal mobile/immobile transport model based on equivalent transformative Caputo formulation}

 \author[label1]{Zhengguang Liu}
  \cortext[cor1]{Corresponding author.}
 \ead{liuzhgsdu@yahoo.com}
 \author[label1]{Xiaoli Li\corref{cor1}}
  \ead{xiaolisdu@163.com}
 \address[label1]{School of Mathematics, Shandong University, Jinan, Shandong 250100, China.}

\begin{abstract}
In this article, we present a new second order finite difference discrete scheme for fractal mobile/immobile transport model based on equivalent transformative Caputo formulation. The new transformative formulation takes the singular kernel away to make the integral calculation more efficient. Furthermore, this definition is also effective where $\alpha$ is a positive integer. Besides, the T-Caputo derivative also helps to increase the convergence rate of the discretization of $\alpha$-order($0<\alpha<1$) Caputo derivative from $O(\tau^{2-\alpha})$ to $O(\tau^{3-\alpha})$, where $\tau$ is the time step. For numerical analysis, a Crank-Nicholson finite difference scheme to solve fractal mobile/immobile transport model is introduced and analyzed. The unconditional stability and a priori estimates of the scheme are given rigorously. Moreover, the applicability and accuracy of the scheme are demonstrated by numerical experiments to support our theoretical analysis.
\end{abstract}

\begin{keyword}


Transformative formulation\sep\ Singular kernel \sep mobile/immobile transport model \sep Unconditional stability \sep Estimates

\MSC[2010] 65M06 \sep 65M12 \sep 65M15 \sep 26A33
\end{keyword}

\end{frontmatter}
 \thispagestyle{empty}
\section{Introductions}
\biboptions{numbers,sort&compress}
In recent years, many problems in physical science, electromagnetism, electrochemistry, diffusion and general transport theory can be solved by the fractional calculus approach, which gives attractive applications as a new modeling tool in a variety of scientific and engineering fields. Roughly speaking, the fractional models can be classified into two principal kinds: space-fractional differential equation and time-fractional one. Numerical methods and theory of solutions of the problems for fractional differential equations have been studied extensively by many researchers which mainly cover finite element methods \cite{zhang2012finite,jiang2011high,li2010finite,zeng2013use}, mixed finite element methods \cite{liu2014new,zhao2015two,liu2014mixed,liu2015h}, finite difference methods \cite{sousa2012second,sousa2014explicit,sousa2009finite,huang2013two}, finite volume (element) methods \cite{cheng2015eulerian,liu2014new1}, (local) discontinuous Galerkin (L)DG methods \cite{wei2014analysis}, spectral methods \cite{lin2007finite,lin2011finite} and so on.

The singular kernel of Caputo fractional derivative causes a lot of difficult problems both in integral calculation and discretization.   To take singular kernel away, Caputo and Fabrizio \cite{caputo2015new} suggest a new definition of fractional derivative by changing the kernel $(t-s)^{-\alpha}$ with the function $exp(-\alpha\frac{t-s}{1-\alpha})$ and $\frac{1}{\Gamma(1-\alpha)}$ with $\frac{M(\alpha)}{1-\alpha}$. The Caputo-Fabrizo derivative can portray substance heterogeneities and configurations with different scales, which noticeably cannot be managing with the renowned local theories. And some related articles have been considered by many authors. Atangana \cite{atangana2016new} introduces the application to nonlinear Fisher¡¯s reaction-diffusion equation based on the new fractional derivative. He \cite{atangana2015extension} also analyzes the extension of the resistance, inductance, capacitance electrical circuit to this fractional derivative without singular kernel. A numerical solution for the model of resistance, inductance, capacitance(RLC) circuit via the fractional derivative without singular kernel is considered by Atangana \cite{atangana2015numerical}. However, we observe that there are many different actions between Caputo-Fabrizio derivative and Caputo derivative. The two definitions are not equivalent and can not transform into each other in any cases.

In this paper, we suggest a new transformative formulation of fractional derivative named T-Caputo forluma, which is equivalent with Caputo fractional derivative in some cases. Furthermore, the two definitions can transform into each other. More importantly, the T-Caputo formula also helps to increase the convergence rate of the discretization of $\alpha$-order($0<\alpha<1$) Caputo derivative from $O(\tau^{2-\alpha})$ to $O(\tau^{3-\alpha})$, where $\tau$ is the time step. For numerical analysis, we present a Crank-Nicholson finite difference scheme to solve fractal mobile/immobile transport model. The unconditional stability and a priori estimates of the scheme are given rigorously. Moreover, the applicability and accuracy of the scheme are demonstrated by numerical experiments to support our theoretical analysis.

A fractal mobile/immobile transport model is a type of second order partial differential equations (PDEs), describing a wide family of problems including heat diffusion and ocean acoustic propagation, in physical or mathematical systems with a time variable, which behave essentially like heat diffusing through a solid \cite{atangana2013numerical}. Significant progress has already been made in the approximation of the time fractional order dispersion equation, see \cite{zhang2009time}. Schumer \cite{schumer2003fractal} firstly developes the fractional-order, mobile/immobile (MIM) model. The time drift term $\partial u/\partial t$ is added to describe the motion time and thus helps to distinguish the status of particles conveniently. This equation is the limiting equation that
governs continuous time random walks with heavy tailed random waiting times. In most cases, it is difficult, or infeasible, to find the analytical solution or good numerical solution of the problems. Numerical solutions or approximate analytical solutions become necessary.
Liu et al. \cite{liu2014rbf} give a radial basis functions(RBFs) meshless approach for modeling a fractal mobile/immobile transport model.
Numerical simulation of the fractional order mobile/immobile advection-dispersion model is consindered by Liu et al. \cite{liu2012numerical}. Furthermore, Zhang and Liu \cite{zhang2013novel} present a novel numerical method for the time variable fractional order mobile--immobile advection--dispersion model. The finite difference schemes are used by Ashyralyev and Cakir \cite{ashyralyev2012numerical} for solving one-dimensional fractional parabolic partial differential equations. They \cite{ashyralyev2013fdm} also give the FDM for fractional parabolic equations with the Neumann condition.

The paper is organized as follows. In Sect.2, we give the definitions and some notations. We introduce a Crank-Nicholson finite difference scheme for a fractal mobile/immobile transport model in Sect.3. Then in Sect.4, we give the analysis of stability and error estimates for the presented method. In Sect.5, some numerical experiments for the second order finite difference discretization are carried out.
\section{Some notations and definitions}
Firstly, we give some definitions which are used in the following analysis.

Let us recall the usual Caputo fractional time derivative of order $\alpha$, given by
\begin{flalign*}
&\leftidx{_0^C}D_t^{\alpha}u(t)=\frac{1}{\Gamma(1-\alpha)}\int_0^tu'(s)(t-s)^{-\alpha}ds,\quad 0<\alpha<1.
\end{flalign*}

Here, we give the following new transformative formulation of fractional derivative.
\begin{defn}\label{def1}
Let $u(t)\in C^2(0,T)$, $\alpha\in(0,1)$, then the new transformative formula of fractional order is defined as:
\begin{equation*}
\leftidx{_0^{TC}}D_t^\alpha u(t)=\frac{1}{\Gamma(2-\alpha)}\int_0^tu''(s)(t-s)^{1-\alpha}ds,\quad 0<\alpha<1.
\end{equation*}
\end{defn}

From the above definition of fractional order transformative formula, we know that the singular kernel $(t-\tau)^{-\alpha}$ in Caputo derivative is replaced with $(t-\tau)^{1-\alpha}$ in new one which does not have singularity for $t=\tau$.

\begin{lem}\label{le1}
Suppose $u(t)\in C^2(0,T)$, $\alpha\in(0,1)$, then we have
\begin{equation*}
\aligned
\leftidx{_0^{TC}}D_t^\alpha u(t)=\leftidx{_0^{C}}D_t^\alpha u(t)-\frac{u'(0)t^{1-\alpha}}{\Gamma(2-\alpha)}.
\endaligned
\end{equation*}
In particular, if the function is such that $u'(0)=0$, then we have
\begin{equation*}
\aligned
\leftidx{_0^{TC}}D_t^\alpha u(t)=\leftidx{_0^{C}}D_t^\alpha u(t).
\endaligned
\end{equation*}
\end{lem}
\textbf{\emph{Proof:}} Noting that
\begin{align*}
\frac{\partial[u'(s)(t-s)^{1-\alpha}]}{\partial s}=u''(s)(t-s)^{1-\alpha}-(1-\alpha)u'(s)(t-s)^{-\alpha}.
\end{align*}
Then it is easy to get
\begin{align*}
u'(s)(t-s)^{-\alpha}=\frac{1}{1-\alpha}\left[u''(s)(t-s)^{1-\alpha}-\frac{\partial[u'(s)(t-s)^{1-\alpha}]}{\partial s}\right].
\end{align*}
Thus the Caputo derivative can be rewritten as
\begin{equation*}
\aligned
\leftidx{_0^{C}}D_t^\alpha u(t)
&=\frac{1}{\Gamma(1-\alpha)}\int_0^tu'(s)(t-s)^{-\alpha}ds\\
&=\frac{1}{\Gamma(2-\alpha)}\int_0^t\left[u''(s)(t-s)^{1-\alpha}-\frac{\partial[u'(s)(t-s)^{1-\alpha}]}{\partial s}\right]ds\\
&=\leftidx{_0^{TC}}D_t^\alpha u(t)-u'(s)(t-s)^{1-\alpha}\left|_0^t\right.\\
&=\leftidx{_0^{TC}}D_t^\alpha u(t)+\frac{u'(0)t^{1-\alpha}}{\Gamma(2-\alpha)}.
\endaligned
\end{equation*}
This completes the proof.

\begin{defn}\label{def2}
Suppose $u(t)\in C^{n+1}(0,T)$, if $n>1$, and $\alpha\in(n-1,n)$, the fractional transformative formulation $\leftidx{_0^{TC}}D_t^\alpha u(t)$ is defined by
\begin{equation*}
\leftidx{_0^{TC}}D_t^\alpha u(t)=\frac{1}{\Gamma(n+1-\alpha)}\int_0^tu^{(n+1)}(s)(t-s)^{n-\alpha}ds,\quad n-1<\alpha<n.
\end{equation*}
\end{defn}

\begin{lem}\label{le2}
Suppose $u(t)\in C^{n+1}(0,T)$, $\alpha\in(n-1,n)$, then we have
\begin{equation*}
\aligned
\leftidx{_0^{TC}}D_t^\alpha u(t)=\leftidx{_0^{C}}D_t^\alpha u(t)-\frac{u^{(n)}(0)t^{n-\alpha}}{\Gamma(n+1-\alpha)}.
\endaligned
\end{equation*}
In particular, if the function is such that $u^{(n)}(0)=0$, then we have
\begin{equation*}
\aligned
\leftidx{_0^{TC}}D_t^\alpha u(t)=\leftidx{_0^{C}}D_t^\alpha u(t),\quad n-1<\alpha<n.
\endaligned
\end{equation*}
\end{lem}
\textbf{\emph{Proof:}} Similarly analysis in the proof of Lemma 1, we have
\begin{align*}
\frac{\partial[u^{(n)}(s)(t-s)^{n-\alpha}]}{\partial s}=u^{(n+1)}(s)(t-s)^{n-\alpha}-(n-\alpha)u^{(n)}(s)(t-s)^{n-1-\alpha}.
\end{align*}
Then it is easy to get
\begin{align*}
u^{(n)}(s)(t-s)^{n-1-\alpha}=\frac{1}{n-\alpha}\left[u^{(n+1)}(s)(t-s)^{n-\alpha}-\frac{\partial[u^{(n)}(s)(t-s)^{n-\alpha}]}{\partial s}\right].
\end{align*}
Thus the $\alpha$-order Caputo derivative can be rewritten as
\begin{equation*}
\aligned
\leftidx{_0^{C}}D_t^\alpha u(t)&=\frac{1}{\Gamma(n-\alpha)}\int_0^tu^{(n)}(s)(t-s)^{n-1-\alpha}ds\\
&=\frac{1}{\Gamma(n+1-\alpha)}\int_0^t\left[u^{(n+1)}(s)(t-s)^{n-\alpha}-\frac{\partial[u^{(n)}(s)(t-s)^{n-\alpha}]}{\partial s}\right]ds\\
&=\leftidx{_0^{TC}}D_t^\alpha u(t)-u^{(n)}(s)(t-s)^{n-\alpha}\left|_0^t\right.\\
&=\leftidx{_0^{TC}}D_t^\alpha u(t)+\frac{u^{(n)}(0)t^{n-\alpha}}{\Gamma(n+1-\alpha)}.
\endaligned
\end{equation*}
This completes the proof.

\begin{lem}\label{le3}
For the new fractional order transformative formulation, $\alpha\in(0,1)$ we have
\begin{equation*}
\aligned
D_t^{(n)}(\leftidx{_0^{TC}}D_t^\alpha u(t))=\leftidx{_0^{C}}D_t^\alpha(D_t^{(n)}u(t)).
\endaligned
\end{equation*}
In particular, if the function is such that $u'(0)=0$, then we have
\begin{equation*}
\aligned
D_t^{(n)}(\leftidx{_0^{TC}}D_t^\alpha u(t))=\leftidx{_0^{TC}}D_t^\alpha(D_t^{(n)}u(t)).
\endaligned
\end{equation*}
\end{lem}
\textbf{\emph{Proof:}} We begin considering $n=1$, then from definition (\ref{def1}) of $\leftidx{_0^{TC}}D_t^\alpha u(t)$, we obtain
\begin{equation*}
\aligned
D_t^{(1)}(\leftidx{_0^{TC}}D_t^\alpha u(t))
&=\frac{d}{dt}\left(\frac{1}{\Gamma(2-\alpha)}\int_0^tu''(s)(t-s)^{1-\alpha}ds\right)\\
&=\frac{1}{\Gamma(2-\alpha)}\left[u''(s)(t-s)^{1-\alpha}\left|_{s=t}\right.+\int_0^t(1-\alpha)u''(s)(t-s)^{-\alpha}ds\right]\\
&=\frac{1}{\Gamma(1-\alpha)}\int_0^tu''(s)(t-s)^{-\alpha}ds\\
&=\leftidx{_0^{C}}D_t^\alpha(D_t^{(1)}u(t)).
\endaligned
\end{equation*}
Particularly, From Lemma \ref{le1}, we know $\leftidx{_0^{TC}}D_t^\alpha u(t)=\leftidx{_0^{C}}D_t^\alpha u(t)$ if $u'(0)=0$. Thus we have
\begin{equation*}
\aligned
D_t^{(1)}(\leftidx{_0^{TC}}D_t^\alpha u(t))=\leftidx{_0^{TC}}D_t^\alpha(D_t^{(1)}u(t)).
\endaligned
\end{equation*}
It is easy to generalize the proof for any $n>1$.

\begin{lem}\label{le3*}
For the new fractional order transformative formulation, if $\alpha=n$, we have
\begin{equation*}
\aligned
{_0^{TC}}D_t^nu(t)=u^{(n)}(t)-u^{(n)}(0).
\endaligned
\end{equation*}
\end{lem}
\textbf{Proof:} From Definition \ref{def2}, we obtain
\begin{equation*}
\aligned
{_0^{TC}}D_t^\alpha u(t)
&=\frac{1}{\Gamma(n+1-\alpha)}\int_0^tu^{(n+1)}(s)(t-s)^{n-\alpha}ds\\.
&=u^{(n)}(s)\left|_{s=0}^t\right.\\
&=u^{(n)}(t)-u^{(n)}(0).
\endaligned
\end{equation*}

From the Lemma \ref{le3*}, we obtain
\begin{equation*}
\aligned
\renewcommand{\arraystretch}{1.5}
{_0^{TC}}D_t^\alpha u(t)
&=\displaystyle\frac{1}{\Gamma(n+1-\alpha)}\int_0^tu^{(n+1)}(s)(t-s)^{n-\alpha}ds\\
&=
  \left\{
   \begin{array}{lr}
\displaystyle{_0^{C}}D_t^\alpha u(t)-\frac{u^{(n)}(0)t^{n-\alpha}}{\Gamma(n+1-\alpha)},&n-1<\alpha<n,\\
u^{(n)}(t)-u^{(n)}(0)&\alpha=n.
\end{array}\right.
\endaligned
\end{equation*}

\begin{figure}[h!b!p!]
\centering
\includegraphics[width=8cm,height=6cm]{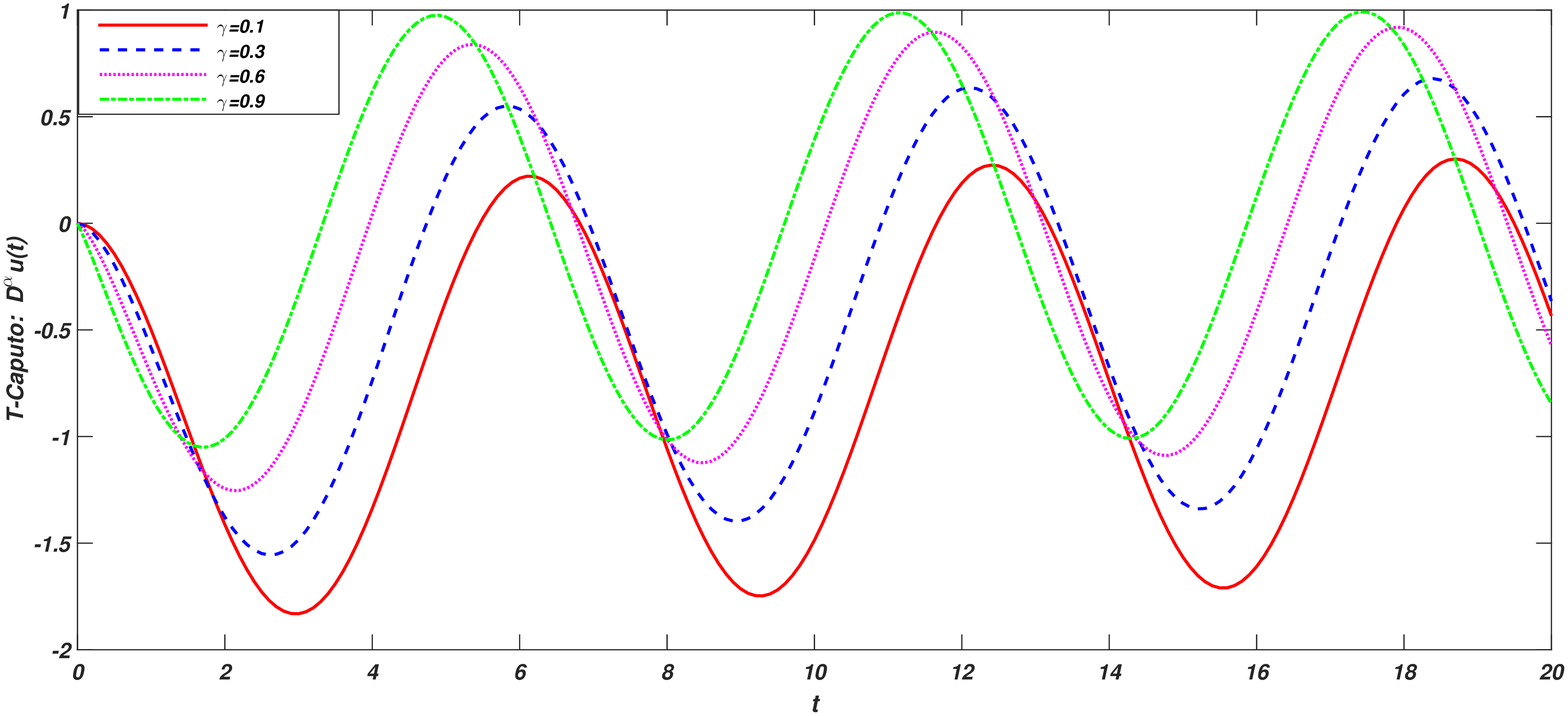}
\includegraphics[width=8cm,height=6cm]{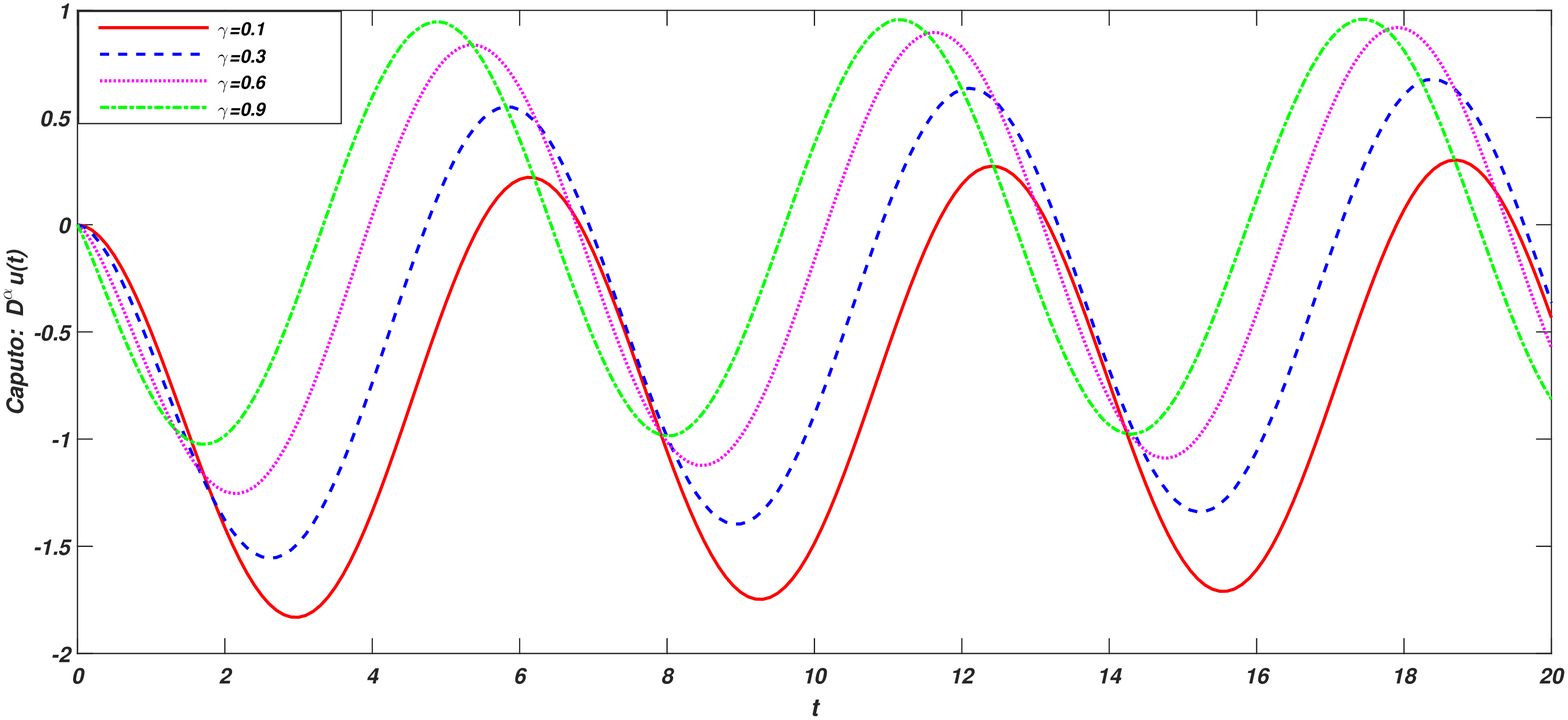}
\caption{Simulation of transformative formulation(left) and Caputo derivative, with $\alpha=$0.1, 0.3, 0.6, 0.9 in the time interval $[0,20]$.}\label{fig1}
\end{figure}

\begin{figure}[h!b!p!]
\centering
\includegraphics[width=8cm,height=6cm]{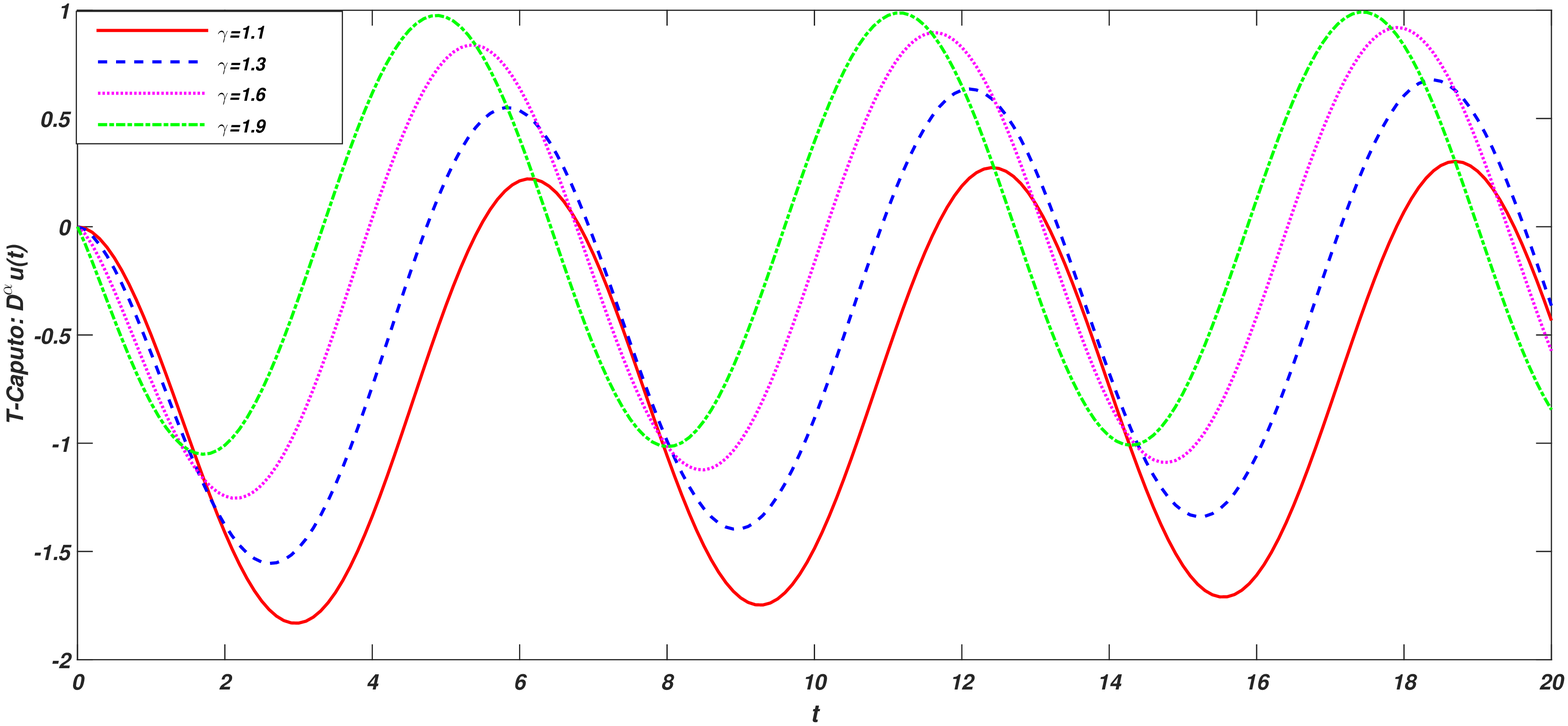}
\includegraphics[width=8cm,height=6cm]{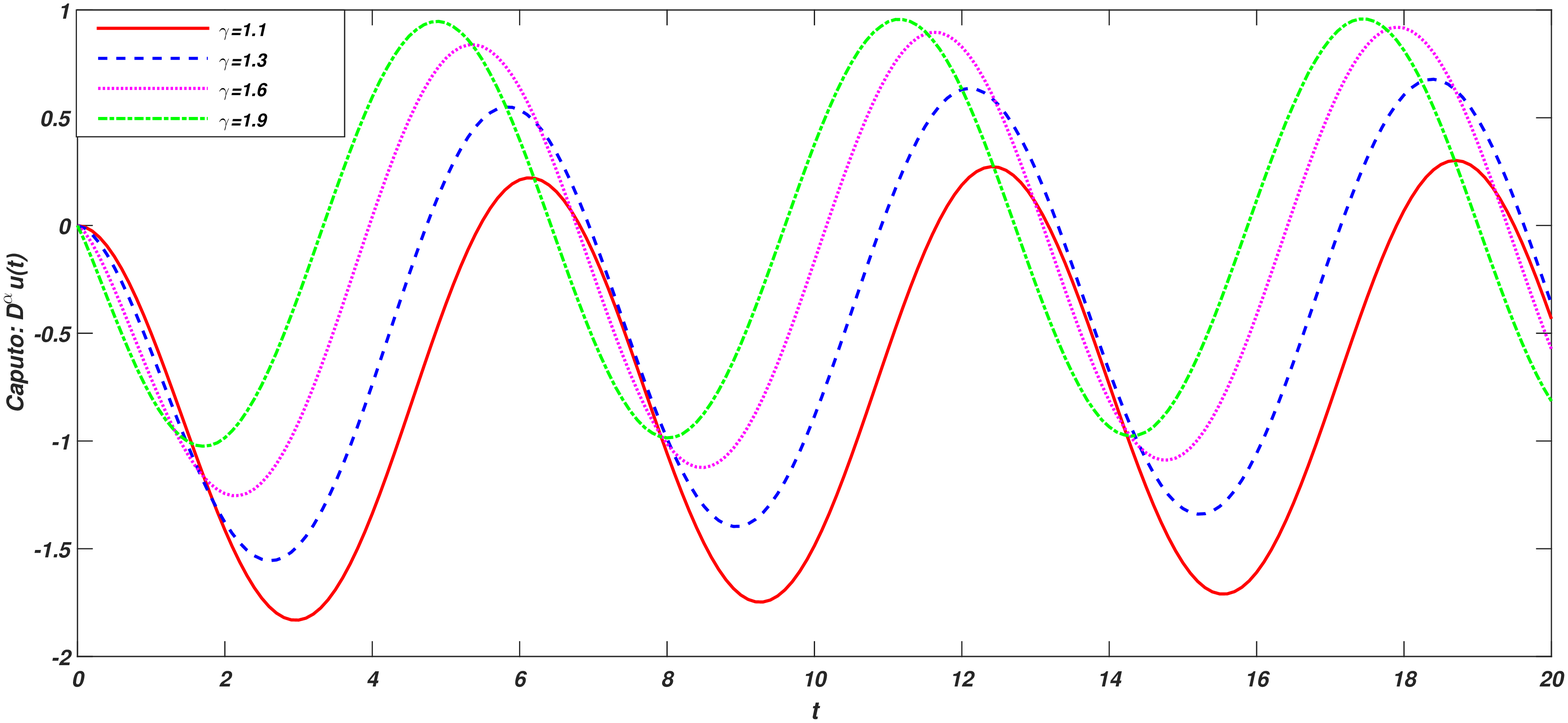}
\caption{Simulation of transformative formulation(left) and Caputo derivative, with $\alpha=$1.1, 1.3, 1.6, 1.9 in the time interval $[0,20]$.}\label{fig2}
\end{figure}

Let us consider, the transformative formulation of a particular function, as $u(t)=cos(t)$ for different $\alpha (0<\alpha<1)$. It is easy to get that $u'(0)=sin(0)=0$. From Figure \ref{fig1}, we observe there are no different actions between transformative formulation and Caputo derivative. We also consider another function as $u(t)=sin(t)$  which has $u''(0)=0$ for different $\alpha (1<\alpha<2)$. From Figure \ref{fig2}, transformative formulation and Caputo derivative have the exact same set of states.

\section{Finite difference scheme for fractal mobile/immobile transport model}
In this section, we introduce the basic ideas for the numerical solution of the fractal mobile/immobile transport model by the second order finite difference scheme.

We consider the following fractal mobile/immobile transport model:
\begin{flalign}\label{e1}
\frac{\partial u(x,t)}{\partial t}+\leftidx{_0^{C}}D_t^\alpha u(x,t)=\frac{\partial^2 u(x,t)}{\partial x^2}+f(x,t),
\end{flalign}
where $(x,t)\in\Omega=[0,L]\times[0,T]$, $0<\alpha<1,$ $f\in C[0,T]$,
with the initial conditions
\begin{flalign}\label{e2}
u(x,0)=\phi(x),\quad 0\leq x\leq L,
\end{flalign}
and boundary conditions
\begin{flalign}\label{e3}
u(0,t)=u(L,t)=0,\quad t>0.
\end{flalign}

Letting $t=0$ in the equation (\ref{e1}), we get
\begin{equation*}
\aligned
u'(0)=\psi(x)=\phi_{xx}(x)+f(x,0).
\endaligned
\end{equation*}
Using Lemma \ref{le1}, the above model can be transformed into the following formulation:
\begin{flalign*}
\renewcommand{\arraystretch}{1.5}
  \left\{
   \begin{array}{l}
\displaystyle\frac{\partial u(x,t)}{\partial t}+\leftidx{_0^{TC}}D_t^\alpha u(x,t)=\frac{\partial^2 u(x,t)}{\partial x^2}+f(x,t)-\frac{\psi(x)t^{1-\alpha}}{\Gamma(2-\alpha)},\quad(x,t)\in\Omega,\\
u(x,0)=\phi(x),\quad 0\leq x\leq L,\\
u(0,t)=u(L,t)=0,\quad t>0,\\
\psi(x)=\phi_{xx}(x)+f(x,0),\quad 0\leq x\leq L.
\end{array}\right.
\end{flalign*}

In order to do discretizations, we define $\Omega_h=\{x_i|x_i=ih,~h=L/M,~0\leq i\leq M\}$  to be a uniform mesh of interval $[0,L]$. Similarly, define $\Omega_\tau=\{t_n,~t_n=i\tau,~\tau=T/N,~0\leq i\leq N\}$ to be a uniform mesh of interval $[0,T]$. The values of the function $u$ at the grid points are denoted $u_j^k=u(x_j,t_k)$. $U_j^k$ is the approximate solution at the point $(x_j,t_k)$. In case, we suppose $V=\{V_i,~0\leq i\leq M, V_0=V_M=0\}$ and $W=\{W_i,~0\leq i\leq N, W_0=W_M=0\}$ are two grid functions on $\Omega_h$. $g=\{g^n,~0\leq n\leq N\}$ is grid functions on $\Omega_\tau$.

For functions $g$, $V$ and $W$, we give some notations, define $L^{2}$ discrete inner products and norms. Define\cite{huang2013two}
\begin{flalign*}
&\delta_tg^n=\frac{g^{n}-g^{n-1}}{\tau},\quad(V,W)=\sum\limits_{i=1}^{M-1}hV_iW_i,\quad \|V\|^2=(V,V).
\end{flalign*}
%

\subsection{The Crank-Nicholson finite difference scheme}
From now on, let $C$ stand for a positive number independent of $\tau$ and $h$, but possibly with different values at different places. We give some lemmas which used in stability analysis and error estimates.

The objective of this section is to consider the Crank-Nicholson finite difference method for equations (\ref{e1}). A discrete approximation to the new transformative formulation $\leftidx{_0^{TC}}D_t^{\alpha}u(x,t)$ at $(x_i,t_{k+\frac{1}{2}})$ can be obtained by the following approximation
\begin{equation}\label{e6}
\aligned
\leftidx{_0^{TC}}D_t^{\alpha}u(x_i,t_{k+\frac{1}{2}})
&=\frac{1}{\Gamma(2-\alpha)}\int_0^{t_{k+\frac{1}{2}}}u^{''}(x_i,s)(t_{k+\frac{1}{2}}-s)^{1-\alpha}ds\\
&=\frac{1}{\Gamma(2-\alpha)}\sum\limits_{j=1}^k\int_{t_{j-\frac{1}{2}}}^{t_{j+\frac{1}{2}}}\left[\frac{u^{'}(x_i,t_{j+\frac{1}{2}})-u^{'}(x_i,t_{j-\frac{1}{2}})}{\tau}+(s-t_{j})u^{(3)}_{t}(x_i,c_j)\right](t_{k+\frac{1}{2}}-s)^{1-\alpha}ds\\
&\quad+\frac{1}{\Gamma(2-\alpha)}\int_{0}^{t_{\frac{1}{2}}}u^{''}(x,s)(t_{k+\frac{1}{2}}-s)^{1-\alpha}ds\\
&=\frac{1}{\Gamma(2-\alpha)}\sum\limits_{j=1}^k\int_{t_{j-\frac{1}{2}}}^{t_{j+\frac{1}{2}}}\left[\frac{u_i^{j+1}-2u_i^j+u_i^{j-1}}{\tau^2}+r^j+(s-t_{j})u^{(3)}_{t}(x_i,c_j)\right](t_{k+\frac{1}{2}}-s)^{1-\alpha}ds\\
&\quad+\frac{1}{\Gamma(2-\alpha)}\int_{-t_{\frac{1}{2}}}^{t_{\frac{1}{2}}}\left[\frac{u_i^{1}-2u_i^0+u_i^{-1}}{\tau^2}+r^0+(s-t_{0})u^{(3)}_{t}(x_i,c_0)\right](t_{k+\frac{1}{2}}-s)^{1-\alpha}ds\\
&\quad-\frac{1}{\Gamma(2-\alpha)}\int_{-t_{\frac{1}{2}}}^{0}u^{''}(x,s)(t_{k+\frac{1}{2}}-s)^{1-\alpha}ds,
\endaligned
\end{equation}
where $c_j\in(x_{j-\frac{1}{2}},x_{j+\frac{1}{2}})$ and for $\xi_1\in(t_{j+\frac{1}{2}},t_{j+1})$, $\xi_2\in(t_{j},t_{j+\frac{1}{2}})$, $\xi_3\in(t_{j-\frac{1}{2}},t_{j})$, $\xi_4\in(t_{j-1},t_{j-\frac{1}{2}})$, $\eta\in(t_{j-\frac{1}{2}},t_{j+\frac{1}{2}})$ and $u(t)\in C^4[0,t_{k+\frac{1}{2}}]$, we have
\begin{equation}\label{e7}
\aligned
r^j=&\frac{1}{24}\tau\left[u^{(3)}_{t}(x_i,t_{j+\frac{1}{2}})-u^{(3)}_{t}(x_i,t_{j-\frac{1}{2}})\right]\\
&+
\frac{1}{256}\tau^2\left[u^{(4)}_{t}(x_i,\xi_1)+u^{(4)}_{t}(x_i,\xi_2)-u^{(4)}_{t}(x_i,\xi_3)-u^{(4)}_{t}(x_i,\xi_4)\right]\\
&=\frac{1}{48}\tau^2u^{(4)}_{t}(x_i,\eta)+O(\tau^2)\\
&=O(\tau^2).
\endaligned
\end{equation}
In particular, for $j=0$, denote $u^{-1}=u^0-\tau u^{'}(x,0)=\phi-\tau\psi$. Using the simple linear interpolant of $u$ at $(-t_1,0)$, so for $s\in(-t_{\frac{1}{2}},0)$, we have $u''(x,s)=0$. It is a suitable method to satisfy the condition $u^{-1}=u^0-\tau u^{'}(x,0)$.

Combining the equation (\ref{e6}) with (\ref{e7}), we obtain
\begin{equation}\label{e9}
\aligned
\leftidx{_0^{TC}}D_t^{\alpha}u(x_i,t_{k+\frac{1}{2}})
&=\frac{\tau^{2-\alpha}}{\Gamma(3-\alpha)}\sum\limits_{j=1}^k\left(\frac{u_i^{j+1}-2u_i^j+u_i^{j-1}}{\tau^2}\right)M_{k-j}
+\frac{\tau^{2-\alpha}}{\Gamma(3-\alpha)}\left(\frac{u_i^{1}-u_i^0}{\tau^2}-\frac{\psi}{\tau}\right)M_k+R_i^{k+\frac{1}{2}}\\
&=\frac{\tau^{1-\alpha}}{\Gamma(3-\alpha)}\left[M_0\delta_tu_i^{k+1}-\sum\limits_{j=1}^k\left(M_{k-j}-M_{k-j+1}\right)\delta_tu_i^j-M_k\psi_i\right]+R_i^{k+\frac{1}{2}},
\endaligned
\end{equation}
where
\begin{equation}\label{e10}
\aligned
M_j=(j+1)^{2-\alpha}-j^{2-\alpha},
\endaligned
\end{equation}
and
\begin{equation}\label{e11}
\aligned
R_i^{k+\frac{1}{2}}&=\frac{1}{\Gamma(2-\alpha)}\sum\limits_{j=0}^k\int_{t_{j-\frac{1}{2}}}^{t_{j+\frac{1}{2}}}\left(r^j+(s-t_{j})u^{(3)}_{t}(x_i,c_j)\right)(t_{k+\frac{1}{2}}-s)^{1-\alpha}ds\\
&=O(\tau^{3-\alpha}).
\endaligned
\end{equation}

We give some Lemmas about $M_j$ that will be used in the following analysis.
\begin{lem}\label{le4}
For the definition $M_j$, $(j=0,1,2,\ldots,N-1)$, we have $M_j>0$ and $M_{j+1}\geq M_j$, $\forall j\leq k$.
\end{lem}
\textbf{Proof:} Observing that $x^{2-\alpha}$ is a monotone increasing function for $0<\alpha<1$, then we have $M_j=(j+1)^{2-\alpha}-j^{2-\alpha}>0$. Next, let $f(x)=(x+1)^{2-\alpha}-x^{2-\alpha}$, we have
\begin{equation*}
\aligned
f^{'}(x)=(2-\alpha)[(x+1)^{1-\alpha}-x^{1-\alpha}]\geq0, \quad\forall x\leq0.
\endaligned
\end{equation*}
Thus we obtain
\begin{equation*}
\aligned
M_{j+1}=f(j+1)\geq f(j)=M_j.
\endaligned
\end{equation*}
This completes the proof.

\begin{lem}\label{le5}
For the definition $M_j=(j+1)^{2-\alpha}-j^{2-\alpha}$, we denote $G_{j+1}=M_{j+1}-M_j$, $(j=0,1,2,\ldots,N-1)$. Then it holds that
\begin{equation*}
G_1\geq G_2\geq \cdots\geq G_N\geq0.
\end{equation*}
\end{lem}
\textbf{Proof:} Firstly, using Lemma (\ref{le4}), it is easy to get $G_j\geq0$. Next, for fixed $0<\alpha<1$, we give the following function
\begin{equation*}
\aligned
f(x)=(x+2)^{2-\alpha}-2(x+1)^{2-\alpha}+x^{2-\alpha},
\endaligned
\end{equation*}
then we have
\begin{equation*}
\aligned
f'(x)=(2-\alpha)[(x+2)^{1-\alpha}-2(x+1)^{1-\alpha}+x^{1-\alpha}].
\endaligned
\end{equation*}
Using Taylor's expansion, we have
\begin{equation*}
\aligned
(x+2)^{1-\alpha}&=(x+1)^{1-\alpha}+(1-\alpha)(x+1)^{-\alpha}-\frac{1}{2!}(1-\alpha)\alpha (x+1)^{-(\alpha+1)}+\frac{1}{3!}(1-\alpha)\alpha(\alpha+1)\xi_1^{-(\alpha+2)},\\
x^{1-\alpha}&=(x+1)^{1-\alpha}-(1-\alpha)(x+1)^{-\alpha}-\frac{1}{2!}(1-\alpha)\alpha (x+1)^{-(\alpha+1)}-\frac{1}{3!}(1-\alpha)\alpha(\alpha+1)\xi_2^{-(\alpha+2)},
\endaligned
\end{equation*}
where $\xi_1\in(x+1,x+2)$ and $\xi_2\in(x,x+1)$.

Thus, we have
\begin{equation*}
\aligned
f'(x)&=(2-\alpha)\left[-(1-\alpha)\alpha (x+1)^{-(\alpha+1)}+\frac{1}{3!}(1-\alpha)\alpha(\alpha+1)\left(\xi_1^{-(\alpha+2)}-\xi_2^{-(\alpha+2)}\right)\right]\\
&\leq-(2-\alpha)(1-\alpha)\alpha(x+1)^{-(\alpha+1)}\\
&\leq0,\quad \forall x\geq0,\quad 0<\alpha<1.
\endaligned
\end{equation*}
It means that $G_j>G_{j+1}$, $\forall j\geq1$. This completes the proof.

The discretization of first order time derivative is stated as:
\begin{equation}\label{e12}
\aligned
&\frac{\partial u(x_i,t_{k+\frac{1}{2}})}{\partial t}=\frac{u_i^{k+1}-u_i^k}{\tau}+O(\tau^2),
\endaligned
\end{equation}
and the second order spatial derivative is stated as:
\begin{equation}\label{e13}
\aligned
&\frac{\partial^2 u(x_i,t_{k+\frac{1}{2}})}{\partial x^2}=\frac{1}{2}\left[\frac{u_{i+1}^{k+1}-2u_i^{k+1}+u_{i-1}^{k+1}}{h^2}+\frac{u_{i+1}^{k}-2u_i^k+u_{i-1}^{k}}{h^2}\right]+O(h^2),
\endaligned
\end{equation}

Combining the equation (\ref{e9}) with equations (\ref{e12})$\sim$(\ref{e13}), we can obtain the following finite difference scheme, $\forall k=0,1,\cdots N-1$,

\begin{equation}\label{e14}
\aligned
&\delta_tU_{i}^{k+1}+\frac{\tau^{1-\alpha}}{\Gamma(3-\alpha)}\left[M_0\delta_tU_i^{k+1}-\sum\limits_{j=1}^k\left(M_{k-j}-M_{k-j+1}\right)\delta_tU_i^j-M_k\psi_i\right]\\
&=\frac{1}{2}\left[\frac{U_{i+1}^{k+1}-2U_i^{k+1}+U_{i-1}^{k+1}}{h^2}+\frac{U_{i+1}^{k}-2U_i^k+U_{i-1}^{k}}{h^2}\right]+f_i^{k+\frac{1}{2}}-\frac{\psi_i[(k+\frac{1}{2})\tau]^{1-\alpha}}{\Gamma(2-\alpha)}.
\endaligned
\end{equation}

Note that $G_{j+1}=M_{j+1}-M_j$, $(j=0,1,2,\ldots,N-1)$, then we have
\begin{equation}\label{e15}
\aligned
&-\sum\limits_{j=1}^k\left(M_{k-j}-M_{k-j+1}\right)\delta_tU_i^j
=\frac{1}{\tau}\left[G_1U_i^k|_{k\geq1}+\sum\limits_{j=1}^{k-1}\left(G_{k-j+1}-G_{k-j}\right)U_i^j-G_kU_i^0|_{k\geq1}\right].
\endaligned
\end{equation}
Let $\beta=\frac{\tau^{1-\alpha}}{\Gamma(3-\alpha)}$, then above scheme (\ref{e14}) can be rewritten as
\begin{equation}\label{e16}
\aligned
&\left[-\frac{\tau}{2h^2}U_{i+1}^{k+1}+\left(\frac{\tau}{h^2}+1+\beta\right)U_i^{k+1}-\frac{\tau}{2h^2}U_{i-1}^{k+1}\right]\\
&=\left[\frac{\tau}{2h^2}U_{i+1}^{k}+\left(-\frac{\tau}{h^2}+1+\beta-G_1|_{k\geq1}\right)U_i^{k}+\frac{\tau}{2h^2}U_{i-1}^{k}+\beta\sum\limits_{j=1}^{k-1}\left(G_{k-j}-G_{k-j+1}\right)U_i^j+G_kU_i^0|_{k\geq1}\right]\\
&\quad+\tau\beta M_k\psi_i+\tau f_i^{k+\frac{1}{2}}-\frac{\tau\psi_i[(k+\frac{1}{2})\tau]^{1-\alpha}}{\Gamma(2-\alpha)}.
\endaligned
\end{equation}
\section{Stability analysis and optimal error estimates}
\subsection{Stability analysis}
We analyze the stability of the difference scheme by a Fourier analysis. Let $\widetilde{U}_i^k$ be the approximate solution of (\ref{e16}), and define
\begin{equation*}
\aligned
\rho_i^k=U_i^k-\widetilde{U}_i^k,\quad 1\leq i\leq M,\quad 0\leq k\leq N.
\endaligned
\end{equation*}
Then we have
\begin{equation}\label{e17}
\aligned
&\left[-\frac{\tau}{2h^2}\rho_{i+1}^{k+1}+\left(\frac{\tau}{h^2}+1+\beta\right)\rho_i^{k+1}-\frac{\tau}{2h^2}\rho_{i-1}^{k+1}\right]\\
&=\left[\frac{\tau}{2h^2}\rho_{i+1}^{k}+\left(-\frac{\tau}{h^2}+1+\beta-\beta G_1|_{k\geq1}\right)\rho_i^{k}+\frac{\tau}{2h^2}\rho_{i-1}^{k}+\beta\sum\limits_{j=1}^{k-1}\left(G_{k-j}-G_{k-j+1}\right)\rho_i^j+\beta G_k\rho_i^0|_{k\geq1}\right].
\endaligned
\end{equation}

As the same definition in \cite{karatay2013new}, we define the grid function
\begin{flalign*}
\rho^k(x)=\left\{
   \begin{array}{ll}
0,&0\leq x\leq x_{\frac{1}{2}},\\
\rho_i^k,&x_{i-\frac{1}{2}}\leq x\leq x_{i+\frac{1}{2}},\quad 1\leq i\leq M-1,\\
0,&x_{M-\frac{1}{2}}\leq x\leq x_{M}.
\end{array}\right.
\end{flalign*}

We can expand $\rho^k(x)$ in a Fourier series
\begin{equation*}
\aligned
\rho^k(x)=\sum\limits_{l=-\infty}^{\infty}d_k(l)e^{\frac{i2\pi lx}{L}},\quad k=1,2,\ldots,N,
\endaligned
\end{equation*}
where discrete Fourier coefficients $d_k(l)$ are
\begin{equation}\label{e18}
\aligned
d_k(l)=\frac{1}{L}\int_0^L\rho^k(\xi)e^{\frac{-i2\pi l\xi}{L}}d\xi.
\endaligned
\end{equation}

Then we have the Parseval equality for the discrete Fourier transform
\begin{equation*}
\aligned
\int_0^L|\rho^k(x)|^2dx=\sum\limits_{l=-\infty}^{\infty}|d_k(l)|^2.
\endaligned
\end{equation*}

Introduce the following norm
\begin{equation*}
\aligned
\|\rho^k\|_2=\left(\sum\limits_{i=1}^{M-1}h|\rho^k_i|^2\right)^{1/2}=\left(\int_0^L|\rho_i^k|^2dx\right)^{1/2}.
\endaligned
\end{equation*}
Then we obtain
\begin{equation*}
\aligned
\|\rho^k\|_2^2=\sum\limits_{l=-\infty}^{\infty}|d_k(l)|^2.
\endaligned
\end{equation*}

Based on the above analysis, we can suppose the solution of equation (\ref{e17}) has the following form $\rho_m^k=d_ke^{imh\gamma}$ where $L=1$ and $\gamma=2\pi l$.
\begin{lem}\label{le6}
Suppose that $d_k(l)$ $(k=1,2,\ldots,N)$ are defined by (\ref{e18}), then for $0<\alpha<1$, we have
\begin{equation*}
|d_k|\leq|d_0|,\quad k=1,2,\ldots,N.
\end{equation*}
\end{lem}
\textbf{Proof:} Substituting $\rho_m^k=d_ke^{imh\gamma}$ into equation (\ref{e17}), we have
\begin{equation}\label{e19}
\aligned
&\left[-\frac{\tau}{2h^2}d_{k+1}e^{i(m+1)h\gamma}+\left(\frac{\tau}{h^2}+1+\beta\right)d_{k+1}e^{imh\gamma}-\frac{\tau}{2h^2}d_{k+1}e^{i(m-1)h\gamma}\right]\\
&=\left[\frac{\tau}{2h^2}d_{k}e^{i(m+1)h\gamma}+\left(-\frac{\tau}{h^2}+1+\beta-\beta G_1|_{k\geq1}\right)d_{k}e^{imh\gamma}+\frac{\tau}{2h^2}d_{k}e^{i(m-1)h\gamma}\right.\\
&\quad\left.+\beta\sum\limits_{j=1}^{k-1}\left(G_{k-j}-G_{k-j+1}\right)d_{j}e^{imh\gamma}+\beta G_kd_{0}e^{imh\gamma}|_{k\geq1}\right].
\endaligned
\end{equation}

By simply calculation, we can get
\begin{equation}\label{e20}
\aligned
\left[-\frac{\tau}{2h^2}(e^{ih\gamma}+e^{-ih\gamma})+\left(\frac{\tau}{h^2}+1+\beta\right)\right]d_{k+1}
=&\left[\frac{\tau}{2h^2}(e^{ih\gamma}+e^{-ih\gamma})+\left(-\frac{\tau}{h^2}+1+\beta-\beta G_1|_{k\geq1}\right)\right]d_k\\
&+\beta\sum\limits_{j=1}^{k-1}\left(G_{k-j}-G_{k-j+1}\right)d_{j}+\beta G_kd_{0}|_{k\geq1}.
\endaligned
\end{equation}
Noting that $e^{ih\gamma}+e^{-ih\gamma}=2\cos(h\gamma)$, thus equation (\ref{e20}) can be rewritten as the following formulation:
\begin{equation}\label{e21}
\aligned
\left[-\frac{\tau}{h^2}\cos(h\gamma)+\left(\frac{\tau}{h^2}+1+\beta\right)\right]d_{k+1}
=&\left[\frac{\tau}{h^2}\cos(h\gamma)+\left(-\frac{\tau}{h^2}+1+\beta-\beta G_1|_{k\geq1}\right)\right]d_k\\
&+\beta\sum\limits_{j=1}^{k-1}\left(G_{k-j}-G_{k-j+1}\right)d_{j}+\beta G_kd_{0}|_{k\geq1}.
\endaligned
\end{equation}

Firstly, letting $k=0$ in equation (\ref{e21}) to obtain
\begin{equation}\label{e22}
\aligned
\displaystyle|d_1|=\left|\frac{\frac{\tau}{h^2}\cos(h\gamma)-\frac{\tau}{h^2}+1+\beta}{-\frac{\tau}{h^2}\cos(h\gamma)+\frac{\tau}{h^2}+1+\beta}\right||d_0|
=\left|\frac{-\frac{1-cos(h\gamma)}{h^2}\tau+1+\beta}{\frac{1-cos(h\gamma)}{h^2}\tau+1+\beta}\right||d_0|\leq|d_0|.
\endaligned
\end{equation}
Now suppose that we have proved that $|d_n|\leq|d_0|$, $n=1,2,\ldots,k$, then using the equation (\ref{e21}), we obtain
\begin{equation}\label{e23}
\aligned
\displaystyle|d_{k+1}|\leq\left(\left|\frac{-\frac{1-cos(h\gamma)}{h^2}\tau+1+\beta-\beta G_1}{\frac{1-cos(h\gamma)}{h^2}\tau+1+\beta}\right||d_k|+\frac{1}{|\frac{1-cos(h\gamma)}{h^2}\tau+1+\beta|}\left[\beta\sum\limits_{j=1}^{k-1}\left(G_{k-j}-G_{k-j+1}\right)|d_j|+\beta G_k|d_0|\right]\right).
\endaligned
\end{equation}

Observing that $G_j\geq0$ and $G_j-G_{j+1}\geq0$ in Lemma \ref{le5}, then we obtain
\begin{equation}\label{e24}
\aligned
\beta\sum\limits_{j=1}^{k-1}\left(G_{k-j}-G_{k-j+1}\right)|d_j|+\beta G_k|d_0|\leq\beta\left(\sum\limits_{j=1}^{k-1}\left(G_{k-j}-G_{k-j+1}\right)+G_k\right)|d_0|=\beta G_1|d_0|
\endaligned
\end{equation}

Combining the equation (\ref{e23}) with equation (\ref{e24}), we can obtain
\begin{equation}\label{e25}
\aligned
\displaystyle|d_{k+1}|\leq\left(\frac{|-\frac{1-cos(h\gamma)}{h^2}\tau+1+\beta-\beta G_1|+\beta G_1}{|\frac{1-cos(h\gamma)}{h^2}\tau+1+\beta|}\right)|d_0|.
\endaligned
\end{equation}
If $-\frac{1-cos(h\gamma)}{h^2}\tau+1+\beta-\beta G_1>0$, then we have
\begin{equation}\label{e26}
\aligned
\displaystyle|d_{k+1}|\leq\left(\frac{-\frac{1-cos(h\gamma)}{h^2}\tau+1+\beta}{|\frac{1-cos(h\gamma)}{h^2}\tau+1+\beta|}\right)|d_0|\leq|d_0|.
\endaligned
\end{equation}
If $-\frac{1-cos(h\gamma)}{h^2}\tau+1+\beta-\beta G_1\leq0$, then we have
\begin{equation}\label{e27}
\aligned
\displaystyle|d_{k+1}|\leq\left(\frac{2\beta G_1+\frac{1-cos(h\gamma)}{h^2}\tau-1-\beta}{|\frac{1-cos(h\gamma)}{h^2}\tau+1+\beta|}\right)|d_0|.
\endaligned
\end{equation}
It means that
\begin{equation}\label{e28}
\aligned
&|d_{k+1}|\leq|d_0|\\
&\Leftrightarrow \frac{2\beta G_1+\frac{1-cos(h\gamma)}{h^2}\tau-1-\beta}{\frac{1-cos(h\gamma)}{h^2}\tau+1+\beta}\leq1\\
&\Leftrightarrow \beta G_1\leq 1+\beta\\
&\Leftrightarrow (2^{2-\alpha}-3)\tau^{1-\alpha}\leq\Gamma(3-\alpha).
\endaligned
\end{equation}

Note that  $(2^{2-\alpha}-3)\tau^{1-\alpha}\leq\Gamma(3-\alpha)$, $\forall\tau\leq1$. It means that $|d_{k+1}|\leq|d_0|$ is unconditionally efficient. By using mathematical induction, we complete the proof.
\begin{thm}\label{le7}
The Crank-Nicholson finite difference scheme defined by (\ref{e16}) is unconditionally stable for $0<\alpha<1$.
\end{thm}
\textbf{Proof:} Suppose that $\widetilde{\textbf{U}}^k$ is the approximate solution of equation (\ref{e16}), Applying Lemma \ref{le6} and Parseval equality, we obtain
\begin{equation}\label{e29}
\aligned
\|\textbf{U}^k-\widetilde{\textbf{U}}^k\|_2^2&=\|\rho^k\|_2^2=\sum\limits_{m=1}^{M-1}h|\rho^k_m|^2=h\sum\limits_{m=1}^{M-1}|d_ke^{imh\gamma}|^2=h\sum\limits_{m=1}^{M-1}|d_k|^2\\
&\leq h\sum\limits_{m=1}^{M-1}|d_0|^2=h\sum\limits_{m=1}^{M-1}|d_0e^{imh\gamma}|^2=\|\rho^0\|_2^2=\|\textbf{U}^0-\widetilde{\textbf{U}}^0\|_2^2,
\endaligned
\end{equation}
which proves that scheme (\ref{e16}) is unconditionally stable.
\subsection{Optimal error estimate}
Combining the equations (\ref{e1}) and (\ref{e9}) with (\ref{e12})$\sim$(\ref{e13}), we obtain
\begin{equation}\label{e30}
\aligned
&\left[-\frac{\tau}{2h^2}u_{i+1}^{k+1}+\left(\frac{\tau}{h^2}+1+\beta\right)u_i^{k+1}-\frac{\tau}{2h^2}u_{i-1}^{k+1}\right]\\
&=\frac{\tau}{2h^2}u_{i+1}^{k}+\left(-\frac{\tau}{h^2}+1+\beta-\beta G_1|_{k\geq1}\right)u_i^{k}+\frac{\tau}{2h^2}u_{i-1}^{k}+\beta\sum\limits_{j=1}^{k-1}\left(G_{k-j}-G_{k-j+1}\right)u_i^j\\
&\quad+\beta G_ku_i^0|_{k\geq1}+\tau\beta M_k\psi_i+\tau f_i^{k+\frac{1}{2}}-\frac{\tau\psi_i[(k+\frac{1}{2})\tau]^{1-\alpha}}{\Gamma(2-\alpha)}+\tau\widetilde{R}_i^{k+\frac{1}{2}},
\endaligned
\end{equation}
where the truncation error at $(x_i,t_{k+\frac{1}{2}})$ is $\widetilde{R}_i^{k+\frac{1}{2}}=O(\tau^2+h^2)$.

Let $\varepsilon_i^k=u_i^k-U_i^k$ be the error at $(x_i,t_{k})$, then subtracting equation (\ref{e16}) from equation (\ref{e30}), we get the error equation as follows
\begin{equation}\label{e31}
\aligned
&\left[-\frac{\tau}{2h^2}\varepsilon_{i+1}^{k+1}+\left(\frac{\tau}{h^2}+1+\beta\right)\varepsilon_i^{k+1}-\frac{\tau}{2h^2}\varepsilon_{i-1}^{k+1}\right]\\
&=\frac{\tau}{2h^2}\varepsilon_{i+1}^{k}+\left(-\frac{\tau}{h^2}+1+\beta-\beta G_1|_{k\geq1}\right)\varepsilon_i^{k}+\frac{\tau}{2h^2}\varepsilon_{i-1}^{k}+\beta\sum\limits_{j=1}^{k-1}\left(G_{k-j}-G_{k-j+1}\right)\varepsilon_i^j\\
&\quad+\beta G_k\varepsilon_i^0|_{k\geq1}+\tau\widetilde{R}_i^{k+\frac{1}{2}},
\endaligned
\end{equation}
Similarly to the stability analysis, we define the grid functions as follows
\begin{flalign*}
\varepsilon^k(x)=\left\{
   \begin{array}{ll}
0,&0\leq x\leq x_{\frac{1}{2}},\\
\varepsilon_i^k,&x_{i-\frac{1}{2}}\leq x\leq x_{i+\frac{1}{2}},\quad 1\leq i\leq M-1,\\
0,&x_{M-\frac{1}{2}}\leq x\leq x_{M},
\end{array}\right.
\end{flalign*}
and
\begin{flalign*}
\widetilde{R}^{k+\frac{1}{2}}(x)=\left\{
   \begin{array}{ll}
0,&0\leq x\leq x_{\frac{1}{2}},\\
\widetilde{R}_i^{k+\frac{1}{2}},&x_{i-\frac{1}{2}}\leq x\leq x_{i+\frac{1}{2}},\quad 1\leq i\leq M-1,\\
0,&x_{M-\frac{1}{2}}\leq x\leq x_{M}.
\end{array}\right.
\end{flalign*}
We can expand $\varepsilon^k(x)$ and $\widetilde{R}^{k+\frac{1}{2}}(x)$ in two Fourier series
\begin{equation*}
\aligned
\varepsilon^k(x)&=\sum\limits_{l=-\infty}^{\infty}\mu_k(l)e^{\frac{2\pi lxi}{L}},\quad k=1,2,\ldots,N,\\
\widetilde{R}^{k+\frac{1}{2}}(x)&=\sum\limits_{l=-\infty}^{\infty}\nu_{k+\frac{1}{2}}(l)e^{\frac{2\pi lxi}{L}},\quad k=0,1,\ldots,N-1,\\
\endaligned
\end{equation*}
where discrete Fourier coefficients $\mu_k(l)$ and $\nu_{k+\frac{1}{2}}(l)$ are
\begin{equation}\label{e32}
\aligned
\mu_k(l)=\frac{1}{L}\int_0^L\varepsilon^k(\xi)e^{\frac{-2\pi l\xi i}{L}}d\xi,\quad \nu_{k+\frac{1}{2}}(l)=\frac{1}{L}\int_0^L\widetilde{R}^{k+\frac{1}{2}}(\xi)e^{\frac{-2\pi l\xi i}{L}}d\xi .
\endaligned
\end{equation}
Then we have the Parseval equality for the discrete Fourier transforms
\begin{equation*}
\aligned
\int_0^L|\varepsilon^k(x)|^2dx=\sum\limits_{l=-\infty}^{\infty}|\mu_k(l)|^2.
\endaligned
\end{equation*}
and
\begin{equation}\label{Pe1}
\aligned
\int_0^L|\widetilde{R}^{k+\frac{1}{2}}(x)|^2dx=\sum\limits_{l=-\infty}^{\infty}|\nu_{k+\frac{1}{2}}(l)|^2.
\endaligned
\end{equation}
Using the boundary conditions, it is easy to obtain $\varepsilon_0^k=\varepsilon_M^k=0$. Thus we define
\begin{equation*}
\aligned
\|\varepsilon^k\|_2=\left(\sum\limits_{i=1}^{M-1}h|\varepsilon^k_i|^2\right)^{1/2}=\left(\int_0^L|\varepsilon_i^k|^2dx\right)^{1/2}.
\endaligned
\end{equation*}
and
\begin{equation*}
\aligned
\|\widetilde{R}^{k+\frac{1}{2}}\|_2=\left(\sum\limits_{i=1}^{M-1}h|\widetilde{R}^{k+\frac{1}{2}}_i|^2\right)^{1/2}=\left(\int_0^L|\widetilde{R}_i^{k+\frac{1}{2}}|^2dx\right)^{1/2}.
\endaligned
\end{equation*}
Without loss of generality, suppose $L=1$, $\gamma=2\pi l$ and
\begin{equation}\label{ee1}
\aligned
\varepsilon_m^k=\mu_ke^{imh\gamma},\quad\widetilde{R}_m^{k+\frac{1}{2}}=\nu_{k+\frac{1}{2}}e^{imh\gamma}.
\endaligned
\end{equation}
Next, Taking notice of the above assumptions (\ref{ee1}), we have
\begin{equation}\label{e33}
\aligned
&\left[-\frac{\tau}{2h^2}\mu_{k+1}e^{i(m+1)h\gamma}+\left(\frac{\tau}{h^2}+1+\beta\right)\mu_{k+1}e^{imh\gamma}-\frac{\tau}{2h^2}\mu_{k+1}e^{i(m-1)h\gamma}\right]\\
&=\frac{\tau}{2h^2}\mu_{k}e^{i(m+1)h\gamma}+\left(-\frac{\tau}{h^2}+1+\beta-\beta G_1|_{k\geq1}\right)\mu_{k}e^{imh\gamma}+\frac{\tau}{2h^2}\mu_{k}e^{i(m-1)h\gamma}\\
&\quad+\beta\sum\limits_{j=1}^{k-1}\left(G_{k-j}-G_{k-j+1}\right)\mu_je^{imh\gamma}+\beta G_k\mu_0|_{k\geq1}e^{imh\gamma}+\tau\nu_{k+\frac{1}{2}}e^{imh\gamma}.
\endaligned
\end{equation}
After simplifications, the equation can be rewritten as
\begin{equation}\label{e34}
\aligned
&\left(\frac{1-\cos(h\gamma)}{h^2}\tau+1+\beta\right)\mu_{k+1}\\
&=\left(-\frac{1-\cos(h\gamma)}{h^2}\tau+1+\beta-\beta G_1|_{k\geq1}\right)\mu_k+\beta\sum\limits_{j=1}^{k-1}\left(G_{k-j}-G_{k-j+1}\right)\mu_j+\beta G_k\mu_0|_{k\geq1}+\tau\nu_{k+\frac{1}{2}}.
\endaligned
\end{equation}
\begin{lem}\label{le8}
Suppose that $\mu_k(l)$ $(k=1,2,\ldots,N)$ and $\nu_{k+\frac{1}{2}}(l)$ $(k=0,1,2,\ldots,N-1)$ are defined by (\ref{e32}), then for $0<\alpha<1$, we have
\begin{equation*}
|\mu_k|\leq C|\nu_{\frac{1}{2}}|,\quad k=1,2,\ldots,N.
\end{equation*}
\end{lem}
\textbf{Proof:} Notice that the error equation satisfies the initial condition $\varepsilon_i^0=0$, $j=0,1,\ldots M$, thus we have $\mu_0=0$.
Firstly, Letting $k=0$, we have
\begin{equation*}
\aligned
&\mu_{1}=\frac{\tau}{\frac{1-\cos(h\gamma)}{h^2}\tau+1+\beta}\nu_{\frac{1}{2}}.
\endaligned
\end{equation*}
It means that $|\mu_1|<|\nu_{\frac{1}{2}}|$.

Now suppose that we have proved that $|\mu_n|\leq C|\mu_0|$, $n=1,2,\ldots,k$, then using the equation (\ref{e34}), we have
\begin{equation}\label{e35}
\aligned
&|\mu_{k+1}|
\leq\left|\frac{-\frac{1-\cos(h\gamma)}{h^2}\tau+1+\beta-\beta G_1}{\frac{1-\cos(h\gamma)}{h^2}\tau+1+\beta}\right||\mu_k|+\frac{\beta\sum\limits_{j=1}^{k-1}\left(G_{k-j}-G_{k-j+1}\right)|\mu_j|+\beta G_k|\mu_0|+\tau|\nu_{k+\frac{1}{2}}|}{\left|\frac{1-\cos(h\gamma)}{h^2}\tau+1+\beta\right|}.
\endaligned
\end{equation}

Similarly to the analysis of equation (\ref{e24}), we obtain
\begin{equation}\label{e36}
\aligned
\beta\sum\limits_{j=1}^{k-1}\left(G_{k-j}-G_{k-j+1}\right)|\mu_j|+\beta G_k|\mu_0|\leq C_1\beta\left(\sum\limits_{j=1}^{k-1}\left(G_{k-j}-G_{k-j+1}\right)+G_k\right)|\nu_{\frac{1}{2}}|=C_1\beta G_1|\nu_{\frac{1}{2}}|.
\endaligned
\end{equation}
Combining the equation (\ref{e35}) with (\ref{e36}), we have
\begin{equation}\label{e37}
\aligned
&|\mu_{k+1}|
\leq\frac{1}{\left|\frac{1-\cos(h\gamma)}{h^2}\tau+1+\beta\right|}\left(C_1\beta G_1+C_2\left|-\frac{1-\cos(h\gamma)}{h^2}\tau+1+\beta-\beta G_1\right|\right)|\nu_{\frac{1}{2}}|+\frac{\tau|\nu_{k+\frac{1}{2}}|}{\left|\frac{1-\cos(h\gamma)}{h^2}\tau+1+\beta\right|}.
\endaligned
\end{equation}
Noting that $\widetilde{R}_{k+\frac{1}{2}}=O(\tau^{2}+h^2)$, $\forall 0\leq k\leq N-1$, and using equation (\ref{Pe1}), we obtain that there is a positive constant $C_{k+\frac{1}{2}}$, such that
\begin{equation*}
\aligned
|\nu_{k+\frac{1}{2}}|\leq C_{k+\frac{1}{2}}|\nu_\frac{1}{2}|, \quad k=0,1,\ldots N-1.
\endaligned
\end{equation*}
Let $C=\max\{C_{\frac{1}{2}},C_{\frac{3}{2}},\ldots,C_{N-\frac{1}{2}}\}$, we have
\begin{equation*}
\aligned
|\nu_{k+\frac{1}{2}}|\leq C|\nu_\frac{1}{2}|, \quad k=0,1,\ldots N-1.
\endaligned
\end{equation*}
Now, let $C=\max(C_1,C_2)$, and if $-\frac{1-\cos(h\gamma)}{h^2}\tau+1+\beta-\beta G_1>0$, then we have
\begin{equation}\label{e38}
\aligned
|\mu_{k+1}|
&\leq C\frac{-\frac{1-\cos(h\gamma)}{h^2}\tau+1+\beta}{\left|\frac{1-\cos(h\gamma)}{h^2}\tau+1+\beta\right|}|\nu_{\frac{1}{2}}|+\frac{C\tau|\nu_{\frac{1}{2}}|}{\left|\frac{1-\cos(h\gamma)}{h^2}\tau+1+\beta\right|}\\
&\leq C|v_{\frac{1}{2}}|.
\endaligned
\end{equation}
If $-\frac{1-\cos(h\gamma)}{h^2}\tau+1+\beta-\beta G_1\leq0$, then we have
\begin{equation}\label{e39}
\aligned
&|\mu_{k+1}|
\leq C\frac{2\beta G_1+\frac{1-\cos(h\gamma)}{h^2}\tau-1-\beta}{\left|\frac{1-\cos(h\gamma)}{h^2}\tau+1+\beta\right|}|\nu_{\frac{1}{2}}|+\frac{C\tau|\nu_{\frac{1}{2}}|}{\left|\frac{1-\cos(h\gamma)}{h^2}\tau+1+\beta\right|}.
\endaligned
\end{equation}
Similarly to the sability analysis, we have
\begin{equation}\label{e40}
\aligned
\frac{2\beta G_1+\frac{1-\cos(h\gamma)}{h^2}\tau-1-\beta}{\left|\frac{1-\cos(h\gamma)}{h^2}\tau+1+\beta\right|}\leq 1,\quad \forall \tau\leq1.
\endaligned
\end{equation}
It means that $|\mu_{k+1}|\leq C|\nu_{\frac{1}{2}}|$. This completes the proof.
\begin{thm}\label{le9}
The Crank-Nicholson finite difference scheme is defined by equation (\ref{e16}) for $0<\alpha<1$, $\textbf{u}^k=(u_1^k,u_2^k,\ldots,u_{M-1}^k)$ and $\textbf{U}^k=(U_1^k,U_2^k,\ldots,U_{M-1}^k)$, then there exists a positive constant $C$ independent of $h$, $k$ and $\tau$ such that.
$$
\|\textbf{u}^k-\textbf{U}^k\|_2\leq C(\tau^2+h^2),\quad \forall 1\leq k\leq N.
$$
\end{thm}
\textbf{Proof:} Applying Lemma \ref{le8} and Parseval equality, we obtain
\begin{equation}\label{e41}
\aligned
\|\textbf{u}^k-\textbf{U}^k\|_2^2&=\|\varepsilon^k\|_2^2=\sum\limits_{m=1}^{M-1}h|\varepsilon^k_m|^2=h\sum\limits_{m=1}^{M-1}|\mu_ke^{imh\gamma}|^2=h\sum\limits_{m=1}^{M-1}|\mu_k|^2\\
&\leq Ch\sum\limits_{m=1}^{M-1}|\nu_{\frac{1}{2}}|^2=Ch\sum\limits_{m=1}^{M-1}|\nu_{\frac{1}{2}}e^{imh\gamma}|^2=C\|\widetilde{R}^{\frac{1}{2}}\|_2^2=C(\tau^2+h^2)^2,
\endaligned
\end{equation}
This completes the proof.

\section{Numerical results}
In this section, some numerical calculations are carried out to test our theoretical results. We consider a numerical example by taking space-time domain $\Omega=[0,1]\times[0,1].$

\textbf{Example 1}:We give the exact solution $u(x,t)=e^tsin(\pi x)$, and for different $\alpha$, we have different $f(x,t)$.
\begin{flalign*}
\renewcommand{\arraystretch}{1.5}
  \left\{
   \begin{array}{l}
\displaystyle\frac{\partial u(x,t)}{\partial t}+\leftidx{_0^{C}}D_t^\alpha u(x,t)=\frac{\partial^2 u(x,t)}{\partial x^2}+f(x,t),\quad(x,t)\in\Omega=[0,1]\times[0,1],\\
u(x,0)=sin(\pi x),\\
u(0,t)=u(1,t)=0.\\
\end{array}\right.
\end{flalign*}

\textbf{Example 2}: The exact solution is $u(x,t)=e^tx^2(1-x)^2$.
\begin{flalign*}
\renewcommand{\arraystretch}{1.5}
  \left\{
   \begin{array}{l}
\displaystyle\frac{\partial u(x,t)}{\partial t}+\leftidx{_0^{C}}D_t^\alpha u(x,t)=\frac{\partial^2 u(x,t)}{\partial x^2}+f(x,t),\quad(x,t)\in\Omega=[0,1]\times[0,1],\\
u(x,0)=x^2(1-x)^2,\\
u(0,t)=u(1,t)=0.\\
\end{array}\right.
\end{flalign*}

\begin{figure}[h!b!p!]
\centering
\includegraphics[width=11cm,height=6cm]{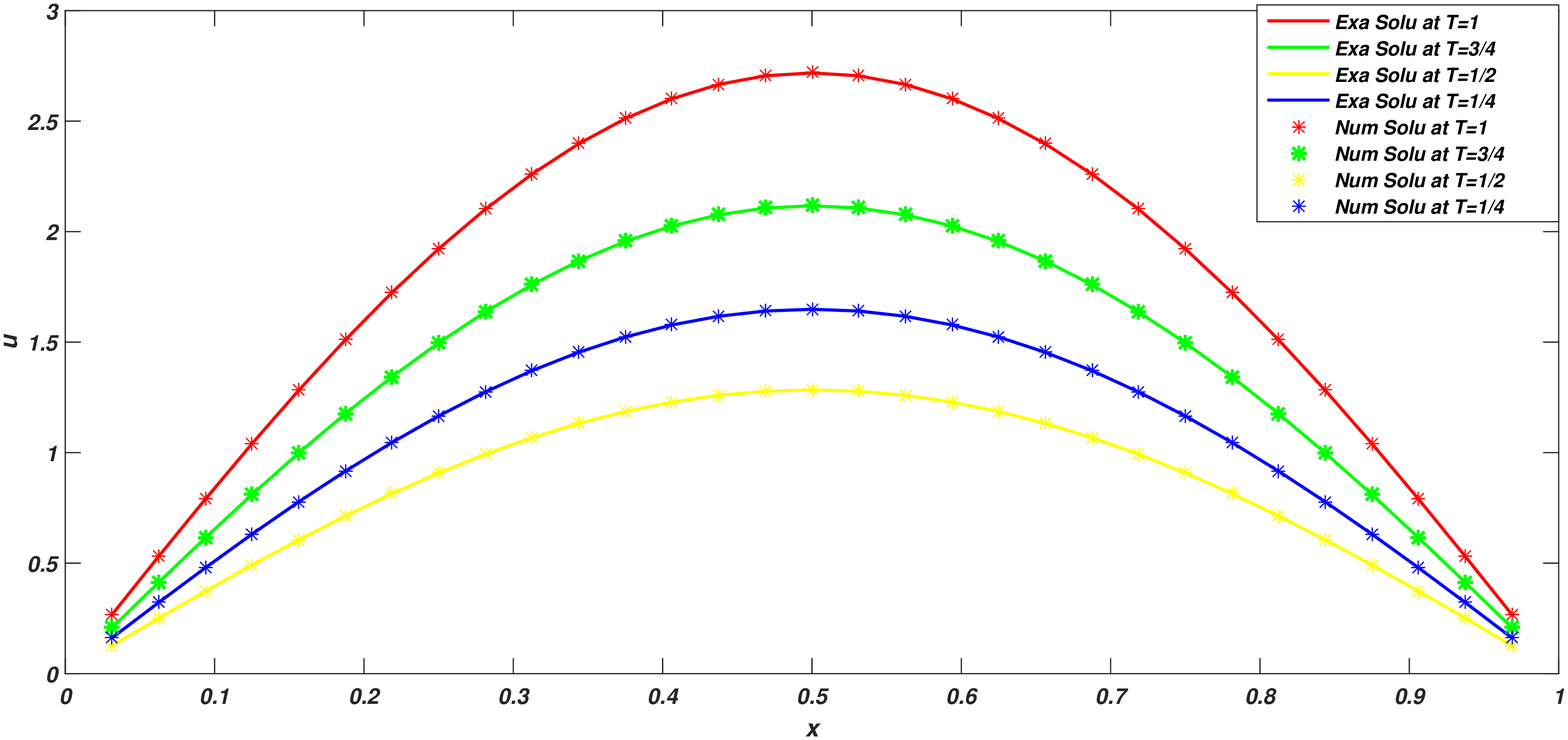}
\caption{Solution behavior of example 1 at $T=1/4$, $T=1/2$, $T=3/4$, $T=1$ with the model parameters $\alpha=1/2$.}\label{fig3}
\end{figure}

\begin{figure}[h!b!p!]
\centering
\includegraphics[width=11cm,height=6cm]{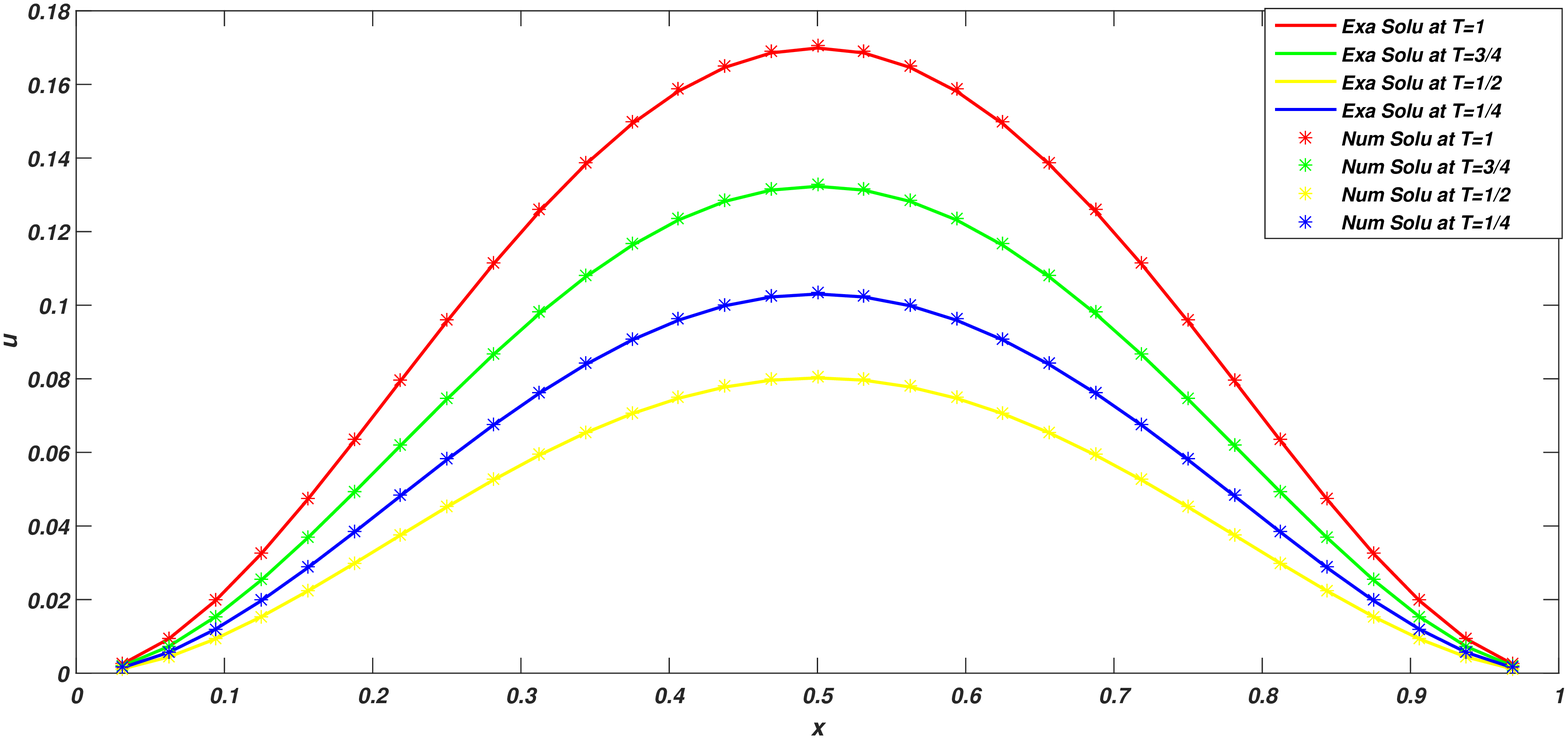}
\caption{Solution behavior of example 2 at $T=1/4$, $T=1/2$, $T=3/4$, $T=1$ with the model parameters $\alpha=1/2$..}\label{fig4}
\end{figure}

Numerical and exact solutions of fractal
mobile/immobile transport model have been depicted in Figure \ref{fig3} (Example 1) and Figure \ref{fig4} (Example 2). Tables \ref{tab1}$\sim$\ref{tab4} show the approximation errors and convergence rates for the second order Crank-Nicholson difference scheme. We take $\tau=\frac{1}{2000}$, a value small enough to check the space errors and convergence rates in Table \ref{tab1} and Table \ref{tab3}. We choose different spatial step sizes to obtain the numerical convergence order in space. In Table \ref{tab2} and Table \ref{tab4}, we take $h=\frac{1}{2000}$, a value small enough such that the spatial discretization errors are negligible as compared with the time errors. we can check that these numerical convergence order almost approaching 2, are consistent with the theoretical analysis.
\begin{table}[h!b!p!]
\small
\renewcommand{\arraystretch}{1.1}
\centering
\caption{\small Errors and convergence rates at the final time $T=1$ of example 1 with different $\alpha$ and take $h=1/2000$.}\label{tab1}
\begin{tabular}{p{0cm}p{0cm}p{1.4cm}p{1cm}p{0.5cm}p{1.4cm}p{1cm}p{0.5cm}p{1.4cm}p{1cm}}
\hline
\multicolumn{1}{l}{\multirow {2}{*}{h}}& \multicolumn{1}{l}{\multirow {2}{*}{$\tau$}}& \multicolumn{2}{l}{$\gamma=0.1$}&&\multicolumn{2}{l}{$\gamma=0.5$}&&\multicolumn{2}{l}{$\gamma=0.9$}\\
\cline{3-4}\cline{6-7} \cline{9-10}
&&\multicolumn{1}{l}{$L_2$ norm error} &\multicolumn{1}{r}{Rate} && \multicolumn{1}{l}{$L_2$ norm error} &\multicolumn{1}{r}{Rate} && \multicolumn{1}{l}{$L_2$ norm error} &\multicolumn{1}{r}{Rate} \\
\hline
\multicolumn{1}{l}{$1/2000$}&\multicolumn{1}{l}{$1/8$}&3.3538e-3&&&3.5156e-3&&&3.3965e-3\\
\multicolumn{1}{l}{$1/2000$}&\multicolumn{1}{l}{$1/16$}&8.8481e-4&\multicolumn{1}{r}{1.9224}&&8.8098e-4&\multicolumn{1}{r}{1.9966}&&8.5199e-4&\multicolumn{1}{r}{1.9951}\\
\multicolumn{1}{l}{$1/2000$}&\multicolumn{1}{l}{$1/32$}&2.2132e-4&\multicolumn{1}{r}{1.9992}&&2.2029e-4&\multicolumn{1}{r}{1.9997}&&2.1334e-4&\multicolumn{1}{r}{1.9977}\\
\multicolumn{1}{l}{$1/2000$}&\multicolumn{1}{l}{$1/64$}&5.5017e-5&\multicolumn{1}{r}{2.0082}&&5.4873e-5&\multicolumn{1}{r}{2.0052}&&5.3220e-5&\multicolumn{1}{r}{2.0031}\\
\multicolumn{1}{l}{$1/2000$}&\multicolumn{1}{l}{$1/128$}&1.3509e-5&\multicolumn{1}{r}{2.0260}&&1.3481e-5&\multicolumn{1}{r}{2.0252}&&1.3089e-5&\multicolumn{1}{r}{2.0236}\\
\hline
\end{tabular}
\end{table}
\begin{table}[h!b!p!]
\small
\renewcommand{\arraystretch}{1.1}
\centering
\caption{\small Errors and convergence rates at the final time $T=1$ of example 1 with different $\alpha$ and take $\tau=1/2000$.}\label{tab2}
\begin{tabular}{p{0cm}p{0cm}p{1.4cm}p{1cm}p{0.5cm}p{1.4cm}p{1cm}p{0.5cm}p{1.4cm}p{1cm}}
\hline
\multicolumn{1}{l}{\multirow {2}{*}{h}}& \multicolumn{1}{l}{\multirow {2}{*}{$\tau$}}& \multicolumn{2}{l}{$\gamma=0.1$}&&\multicolumn{2}{l}{$\gamma=0.5$}&&\multicolumn{2}{l}{$\gamma=0.9$}\\
\cline{3-4}\cline{6-7} \cline{9-10}
&&\multicolumn{1}{l}{$L_2$ norm error} &\multicolumn{1}{r}{Rate} && \multicolumn{1}{l}{$L_2$ norm error} &\multicolumn{1}{r}{Rate} && \multicolumn{1}{l}{$L_2$ norm error} &\multicolumn{1}{r}{Rate} \\
\hline
\multicolumn{1}{l}{$1/8$}&\multicolumn{1}{l}{$1/2000$}&2.0603e-2&&&2.0487e-2&&&2.0542e-2\\
\multicolumn{1}{l}{$1/16$}&\multicolumn{1}{l}{$1/2000$}&5.1295e-3&\multicolumn{1}{r}{2.0060}&&5.1010e-3&\multicolumn{1}{r}{2.0059}&&5.1149e-3&\multicolumn{1}{r}{2.0058}\\
\multicolumn{1}{l}{$1/32$}&\multicolumn{1}{l}{$1/2000$}&1.2810e-3&\multicolumn{1}{r}{2.0015}&&1.2739e-3&\multicolumn{1}{r}{2.0015}&&1.2774e-3&\multicolumn{1}{r}{2.0015}\\
\multicolumn{1}{l}{$1/64$}&\multicolumn{1}{l}{$1/2000$}&3.2012e-4&\multicolumn{1}{r}{2.0006}&&3.1836e-4&\multicolumn{1}{r}{2.0005}&&3.1923e-4&\multicolumn{1}{r}{2.0005}\\
\multicolumn{1}{l}{$1/128$}&\multicolumn{1}{l}{$1/2000$}&7.9984e-5&\multicolumn{1}{r}{2.0008}&&7.9542e-5&\multicolumn{1}{r}{2.0009}&&7.9761e-5&\multicolumn{1}{r}{2.0008}\\
\hline
\end{tabular}
\end{table}

\begin{table}[h!b!p!]
\small
\renewcommand{\arraystretch}{1.1}
\centering
\caption{\small Errors and convergence rates at the final time $T=1$ of example 2 with different $\alpha$ and take $h=1/2000$..}\label{tab3}
\begin{tabular}{p{0cm}p{0cm}p{1.4cm}p{1cm}p{0.5cm}p{1.4cm}p{1cm}p{0.5cm}p{1.4cm}p{1cm}}
\hline
\multicolumn{1}{l}{\multirow {2}{*}{h}}& \multicolumn{1}{l}{\multirow {2}{*}{$\tau$}}& \multicolumn{2}{l}{$\gamma=0.1$}&&\multicolumn{2}{l}{$\gamma=0.5$}&&\multicolumn{2}{l}{$\gamma=0.9$}\\
\cline{3-4}\cline{6-7} \cline{9-10}
&&\multicolumn{1}{l}{$L_2$ norm error} &\multicolumn{1}{r}{Rate} && \multicolumn{1}{l}{$L_2$ norm error} &\multicolumn{1}{r}{Rate} && \multicolumn{1}{l}{$L_2$ norm error} &\multicolumn{1}{r}{Rate} \\
\hline
\multicolumn{1}{l}{$1/2000$}&\multicolumn{1}{l}{$1/8$}&1.9924e-4&&&1.9816e-4&&&1.9166e-4\\
\multicolumn{1}{l}{$1/2000$}&\multicolumn{1}{l}{$1/16$}&4.9831e-5&\multicolumn{1}{r}{1.9994}&&4.9621e-5&\multicolumn{1}{r}{1.9976}&&4.8024e-5&\multicolumn{1}{r}{1.9967}\\
\multicolumn{1}{l}{$1/2000$}&\multicolumn{1}{l}{$1/32$}&1.2386e-5&\multicolumn{1}{r}{2.0083}&&1.2346e-5&\multicolumn{1}{r}{2.0069}&&1.1963e-5&\multicolumn{1}{r}{2.0052}\\
\multicolumn{1}{l}{$1/2000$}&\multicolumn{1}{l}{$1/64$}&3.0214e-6&\multicolumn{1}{r}{2.0354}&&3.0140e-6&\multicolumn{1}{r}{2.0343}&&2.9227e-6&\multicolumn{1}{r}{2.0332}\\
\multicolumn{1}{l}{$1/2000$}&\multicolumn{1}{l}{$1/128$}&6.8010e-7&\multicolumn{1}{r}{2.1514}&&6.7944e-7&\multicolumn{1}{r}{2.1493}&&6.5733e-7&\multicolumn{1}{r}{2.1526}\\
\hline
\end{tabular}
\end{table}
\begin{table}[h!b!p!]
\small
\renewcommand{\arraystretch}{1.1}
\centering
\caption{\small Errors and convergence rates at the final time $T=1$ of example 2 with different $\alpha$ and take $\tau=1/2000$.}\label{tab4}
\begin{tabular}{p{0cm}p{0cm}p{1.4cm}p{1cm}p{0.5cm}p{1.4cm}p{1cm}p{0.5cm}p{1.4cm}p{1cm}}
\hline
\multicolumn{1}{l}{\multirow {2}{*}{h}}& \multicolumn{1}{l}{\multirow {2}{*}{$\tau$}}& \multicolumn{2}{l}{$\gamma=0.1$}&&\multicolumn{2}{l}{$\gamma=0.5$}&&\multicolumn{2}{l}{$\gamma=0.9$}\\
\cline{3-4}\cline{6-7} \cline{9-10}
&&\multicolumn{1}{l}{$L_2$ norm error} &\multicolumn{1}{r}{Rate} && \multicolumn{1}{l}{$L_2$ norm error} &\multicolumn{1}{r}{Rate} && \multicolumn{1}{l}{$L_2$ norm error} &\multicolumn{1}{r}{Rate} \\
\hline
\multicolumn{1}{l}{$1/8$}&\multicolumn{1}{l}{$1/2000$}&6.4192e-3&&&6.3831e-3&&&6.4003e-3\\
\multicolumn{1}{l}{$1/16$}&\multicolumn{1}{l}{$1/2000$}&1.6076e-3&\multicolumn{1}{r}{1.9975}&&1.5987e-3&\multicolumn{1}{r}{1.9974}&&1.6031e-3&\multicolumn{1}{r}{1.9973}\\
\multicolumn{1}{l}{$1/32$}&\multicolumn{1}{l}{$1/2000$}&4.0208e-4&\multicolumn{1}{r}{1.9994}&&3.9986e-4&\multicolumn{1}{r}{1.9993}&&4.0095e-4&\multicolumn{1}{r}{1.9994}\\
\multicolumn{1}{l}{$1/64$}&\multicolumn{1}{l}{$1/2000$}&1.0052e-4&\multicolumn{1}{r}{2.0000}&&9.9758e-5&\multicolumn{1}{r}{2.0030}&&1.0024e-4&\multicolumn{1}{r}{2.0000}\\
\multicolumn{1}{l}{$1/128$}&\multicolumn{1}{l}{$1/2000$}&2.5130e-5&\multicolumn{1}{r}{2.0000}&&2.4992e-5&\multicolumn{1}{r}{1.9970}&&2.5060e-5&\multicolumn{1}{r}{2.0000}\\
\hline
\end{tabular}
\end{table}

\section{Conclusion}
 In this article, we define a novel transformative Caputo fractional derivative which is equivalent with Caputo fractional derivative. This new transformative Caputo derivative takes the singular kernel away to make the integral calculation more efficient. Furthermore, the transformative formulation also helps to increase the convergence rate of the discretization of $\alpha$-order($0<\alpha<1$) Caputo derivative from $O(\tau^{2-\alpha})$ to $O(\tau^{3-\alpha})$, where $\tau$ is the time step. We prove some lemmas and give a Crank-Nicholson finite difference scheme for fractal
mobile/immobile transport model. By using transformative formulation, second-order error estimates in both of temporal and spatial mesh-size in descrete $L^\infty(L^2)$ errors are established for the Crank-Nicholson finite difference scheme.
\bibliographystyle{model1a-num-names}
\bibliography{rough}
\end{document}